\documentclass[11pt]{amsart}
\usepackage[usenames]{color}
\usepackage{epsf}
\usepackage{fullpage}
\usepackage{hyperref}
\usepackage{tikz}
\usetikzlibrary{arrows,calc,positioning,through,intersections,backgrounds,matrix,fit,shapes,decorations}
\newcommand{\avoidbreak}{\postdisplaypenalty=100}
\newtheorem{theorem}{Theorem}[section]
\newtheorem{defn}[theorem]{Definition}

\newtheorem{lemma}[theorem]{Lemma}

\newtheorem{claim}[theorem]{Claim}

\title{Tight co-degree condition for the existence of loose Hamilton cycles in 3-graphs}

\newcommand{\DownEdge}[1]{
  \def\first{#1}
  \pgfmathtruncatemacro\second{\first + 1};
  \pgfmathtruncatemacro\third{\first + 2};
  \node [vertex,below right=8pt of v\first] (v\second) {};
  \node [vertex,below right=8pt of v\second] (v\third) {};
  \node [rotate fit=-45,draw=black,inner sep=5pt,shape=rounded
  rectangle,fit=(v\first) (v\second) (v\third)] {};
}

\newcommand{\DownEdgeWhite}[1]{
  \def\first{#1}
  \pgfmathtruncatemacro\second{\first + 1};
  \pgfmathtruncatemacro\third{\first + 2};
  \node [vertex,below right=8pt of v\first] (v\second) {};
  \node [vertex,fill=white,draw=black,below right=8pt of v\second] (v\third) {};
  \node [rotate fit=-45,draw=black,inner sep=5pt,shape=rounded
  rectangle,fit=(v\first) (v\second) (v\third)] {};
}

\newcommand{\UpEdge}[1]{
  \def\first{#1}
  \pgfmathtruncatemacro\second{\first + 1};
  \pgfmathtruncatemacro\third{\first + 2};
  \node [vertex,above right=8pt of v\first] (v\second) {};
  \node [vertex,above right=8pt of v\second] (v\third) {};
  \node [rotate fit=45,draw=black,inner sep=5pt,shape=rounded
  rectangle,fit=(v\first) (v\second) (v\third)] {};
}

\newcommand{\UpEdgeWhite}[1]{
  \def\first{#1}
  \pgfmathtruncatemacro\second{\first + 1};
  \pgfmathtruncatemacro\third{\first + 2};
  \node [vertex,above right=8pt of v\first] (v\second) {};
  \node [vertex,draw=black,fill=white,above right=8pt of v\second] (v\third) {};
  \node [rotate fit=45,draw=black,inner sep=5pt,shape=rounded
  rectangle,fit=(v\first) (v\second) (v\third)] {};
}

\newcommand{\TwoPath}[1]{
  \def\first{#1}
  \pgfmathtruncatemacro\second{\first + 1};
  \pgfmathtruncatemacro\third{\first + 2};
  \pgfmathtruncatemacro\fourth{\first + 3};
  \pgfmathtruncatemacro\fifth{\first + 4};
  \node [vertex,right=8pt of v\first] (v\second) {};
  \node [vertex,right=8pt of v\second] (v\third) {};
  \node [vertex,right=8pt of v\third] (v\fourth) {};
  \node [vertex,right=8pt of v\fourth] (v\fifth) {};
  \node [draw=black,inner sep=0pt,shape=ellipse,fit=(v\first) (v\second)
  (v\third)] {};
  \node [draw=black,inner sep=0pt,shape=ellipse,fit=(v\third) (v\fourth)
  (v\fifth)] {};
}

\newcommand{\UpDownEdge}[1]{
  \def\first{#1}
  \pgfmathtruncatemacro\second{\first + 1};
  \pgfmathtruncatemacro\third{\first + 2};

  \node [vertex,below right=12pt of v\first] (v\second) {};
  \node [vertex,above right=12pt of v\second] (v\third) {};
  \node [draw=black,inner sep=1pt, 
  shape=ellipse,inner sep=1pt,fit=(v\first) (v\second) (v\third)] {};
}

\begin{document}

\author{Andrzej Czygrinow$^{*}$}

\thanks{$^{*}$School of Mathematical Sciences and Statistics, Arizona
State University, Tempe, AZ 85287, USA. E-mail address: aczygri@asu.edu.
Research of this author is partially supported by NSA grant H98230-13-1-0211.}

\author{Theodore Molla$^{\dagger}$}

\thanks{$^{\dagger}$School of Mathematical Sciences and Statistics, Arizona
State University, Tempe, AZ 85287, USA. E-mail address: tmolla@asu.edu.
Research of this author is partially supported by NSA grant H98230-12-1-0212.}

\begin{abstract}
In 2006, K\"{u}hn and Osthus showed that if a 3-graph $H$ on $n$ vertices has
minimum co-degree at least $(1/4 +o(1))n$ and $n$ is even then 
$H$ has a loose Hamilton cycle. In this paper, we prove that the minimum co-degree of $n/4$ suffices. The result is tight.
\end{abstract}

\maketitle

\section{Introduction}
A $k$-uniform hypergraph $H$ is a pair $(V,E)$ such that $V$ is a finite set
of vertices and the edge set $E$ is a subset of 
$\binom{V}{k}$. Often, we will identify $H$
with its edges and if needed we will use $V(H)$ and $E(H)$ to denote the vertex
set and the edge set of $H$ respectively. 
In this paper we will only use 3-uniform
hypergraphs (3-graphs) and 2-uniform hypergraphs (graphs). We say that a
3-graph $H$ is a \emph{(loose) path} 
if its vertices can be ordered as $v_1, v_2,
\dots,v_{2m+1}$ so that $H= \{\{v_{2i+1}, v_{2i+2}, v_{2i+3}\}: i=0, \dots,
m-1\}$ with \emph{endpoints} $\{v_1, v_{2m+1}\}$. 
$H$ is called a \emph{(loose) cycle} if its vertices can be ordered as
$v_1, v_2, \dots,v_{2m}$  so that $$H=  \{\{v_{2i+1}, v_{2i+2}, v_{2i+3}\}:
i=0, \dots, m-2\} \cup \{\{v_1, v_{2m}, v_{2m-1}\}\}.$$  We say that a 3-graph
$H$ contains a loose Hamilton cycle if $H$ has a loose cycle containing
$V(H)$.  Clearly, for $H$  to have a Hamilton cycle the order of $H$ must be
even and from now on we will always assume that this is the case. Note that in
addition to loose cycles different concepts of Hamilton cycles can be
considered.  Most notably, tight cycles in which every block of three
consecutive vertices is a hyperedge have been intensively studied (see
\cite{RRS3}).

There is a substantial body of work on Hamilton problems in graphs. One direction that has been a source of many important  theorems asks for sufficient minimum degree conditions. Among many theorems of that flavor, Dirac's theorem
is one of the cornerstones of graph theory. In the realm of 3-graphs different
notions of neighborhood and hence minimum degree can be considered. 
For a $3$-graph $H = (V, E)$ and distinct $v, v' \in V$ define 
$N(v) := \{ \{u, u'\} \in \binom{V}{2} : \{v, u, u'\} \in H \}$
and 
$N(v, v') := \{ u \in V : \{v, v', u\} \in H \}$.
In this paper, we will be interested in 
the \emph{minimum co-degree} of $H$ which is 
$\min_{\{v,v'\} \in \binom{V}{2}} \{ |N(v, v')| \}$.
For 3-graphs and loose cycles, K\"{u}hn and Osthus proved the following theorem in \cite{KO}.
\begin{theorem}\label{KO-thm}
For each $\epsilon> 0$ there is an integer $n_0$ such that every 3-graph $H$ on $n\geq n_0$ vertices where $n$ is even and with minimum co-degree at least $n/4 + \epsilon n$ contains a loose Hamilton cycle.
\end{theorem} 
In  the case of tight cycles the minimum co-degree must be at least $(1/2
+o(1))n$ as showed by R\"{o}dl, Ruci\'{n}ski, and Szemer\'{e}di in \cite{RRS3}.
 
 In this paper we get rid of the $\epsilon n$ factor in Theorem~\ref{KO-thm} and show that $n/4$ suffices.
 \begin{theorem}\label{main-thm}
There is an integer $n_0$ such that every 3-graph $H$ on $n\geq n_0$ vertices where $n$ is even and with minimum co-degree at least $n/4$ contains a loose Hamilton cycle.
\end{theorem} 
It is not difficult to see that Theorem~\ref{main-thm} is tight. Indeed, for
even $n$ not divisible by four, consider $H^*$ with $V(H^*)= U \cup W$ such
that $|U|= \frac{n-2}{4}$, $|W|=\frac{3n+2}{4}$ and $H^*$ containing all edges
that intersect $U$. It is not hard to see that $H^*$ has no loose Hamilton cycle. 

Our proof of Theorem~\ref{main-thm} uses the so-called stability approach. In
this method we consider two main cases: the extremal and the non-extremal. In
the non-extremal case we show that if $H$ is ``far'' from $H^*$ and satisfies
the assumptions of the theorem, then $H$ has a loose Hamilton cycle. In the
extremal case we handle the case when $H$ looks like $H^*$. More formally, the
following notion will be used to split the proof into two cases.
\begin{defn}
A 3-graph $H$ is called $\beta$-extremal if it is possible to  find a set
$W\subset V(H)$ such that $|W| = \left\lceil \frac{3n}{4} \right\rceil$ and the 3-graph induced by $W$, $H[W]$, has at most $\beta |V|^3$ edges.
\end{defn} 

The rest of the paper is organized as follows. After a short description
of notation and terminology, we prove the non-extremal
case.  In the last section we establish the extremal case.

\begin{subsection}{Notation and terminology}
In addition to the already introduced terminology we need 
a few concepts. 
For a $k$-graph $H = (V, E)$ and
$V_1, \dotsc, V_k \subseteq V$ define
$E_H(V_1, \dotsc, V_k) := 
\{\{v_1, \dotsc, v_k\} \in H : v_i \in V_i \text{ for $i = 1, \dotsc, k$} \}$
and $e_H(V_1, \dotsc , V_k) := |E_H(V_1, \dotsc, V_k)|$. 
When it is clear from context we may drop the subscript $H$.
For pairwise disjoint $V_1, \dotsc, V_k$ 
let $H[V_1, \dotsc, V_k] := (V_1 \cup \dotsm \cup V_k, E_H(V_1, \dotsc, V_k))$
be the \emph{$k$-partite graph induced by $\{V_1, \dotsc, V_k\}$}.

The following definitions are only used in the proof of Lemma~\ref{main-extremal}.
Let $G$ be a bipartite graph with partite sets $V_1$ and $V_2$.  For
$X_1 \subseteq V_1$ and $X_2 \subseteq V_2$ both non-empty,
define $d(X_1, X_2) := \frac{e_G(X_1, X_2)}{|X_1||X_2|}.$
For constants $0 < \epsilon < 1$ and $0 \le d \le 1$,
we say that $G$ is \emph{$(d, \epsilon)$-regular} if 
$$(1 - \epsilon)d \le d(X_1, X_2) \le (1 + \epsilon)d$$
whenever $|X_i| \ge \epsilon |V_i|$ for $i = 1, 2$.
We say that $G$ is \emph{$(d, \epsilon)$-superregular} if
$G$ is $(d, \epsilon)$-regular and
$(1 - \epsilon)d|V_i| \leq \deg(v, V_i) \leq (1+\epsilon)d|V_i|$ for
every $v \in V_{3-i}$ and $i = 1,2$. The preceding two definitions
match the corresponding definitions in \cite{KO-match}.
Also, we only use these definitions with $d = 1$, so the upper bounds
in the definitions are trivially satisfied.
We will assume throughout that $n$ is sufficiently large.
\end{subsection}
 
\begin{section}{Non-extremal case}
In this section, we prove the non-extremal case. 
\begin{theorem}\label{non-extremal}
For every $\beta>0$ there is an integer $n_0$ such that every 3-graph 
$H = (V, E)$ on $n\geq n_0$ vertices where $n$ is even which is not $\beta$-extremal and with minimum co-degree at least $n/4$ contains a loose Hamilton cycle.
\end{theorem}
\noindent{\bf Proof.} To prove Theorem~\ref{non-extremal} we will use the
connecting-absorbing-reservoir method as in \cite{RRS3}.
\begin{defn}
Let $v, v'$ be two vertices in $V$. A set $U=\{u_1, \dots, u_5\}$ is said to
absorb $v,v'$ if there is a path in $H[U]$ on five vertices and a path in $H[U
  \cup \{v, v'\}]$ on seven vertices both with endpoints $u_1$ and $u_5$.
\end{defn}
\begin{claim}\label{lem1}
Let $v, v'$ be two distinct vertices in $V$. Then there are at least $\binom{n}{5}/500$ sets in $H$ that absorb $v,v'$.
\end{claim}
\noindent{\bf Proof.} We can select $\{u_1, \dots, u_5\} = U$ as follows. Let
$u_1 \in V \setminus \{v, v'\}$ be an arbitrary vertex. Select $u_2 \in N(u_1,
v) \setminus \{v'\}$, $u_3 \in N(u_1, u_2) \setminus \{v,v'\}$, 
$u_4 \in N(v', u_2) \setminus \{v, u_1, u_3\}$ and $u_5 \in N(u_3, u_4) \setminus \{v, v', u_1, u_2\}$. Then $\{u_1, u_2, u_3\}, \{u_3, u_4, u_5\} \in H$ and $\{u_1, v, u_2\}, \{u_2, v', u_4\}, \{u_4, u_3, u_5\}\in H$
(see Figure~\ref{fig:absorber}).
There are at least $(n-2)(n/4 -1)(n/4-2)(n/4-3)(n/4-4)/5! \geq
\binom{n}{5}/500$ choices for $U$. \qed
\begin{figure}
  \centering
  \scalebox{0.7} {
  \begin{tikzpicture}[
    vertex/.style={circle,draw=black,minimum size=24pt},
    edge/.style={double distance=30pt,line cap=round,opacity=0.25}]
    \node [vertex,] (u1) {$u_1$};
    \node [vertex,above right=of u1] (u2) {$u_2$};
    \node [vertex,above right=of u2] (u3) {$u_3$};
    \node [vertex,below right=of u3] (u4) {$u_4$};
    \node [vertex,below right=of u4] (u5) {$u_5$};

    \node [vertex,above=of u1] (v1) {$v$};
    \node [vertex,right=of u2] (v2) {$v'$};

    \begin{pgfonlayer}{background}
      \node [rotate fit=45,ultra thick,draw=black,inner sep=5pt,shape=rounded rectangle,fit=(u1) (u2) (u3)] {};
      \node [rotate fit=-45,ultra thick,draw=black,inner sep=5pt,shape=rounded rectangle,fit=(u3) (u4) (u5)] {};
      \node [draw=black,ultra thick,inner sep=3pt,shape=rounded rectangle,fit=(u2) (v2) (u4)] {};
      \node [draw=black,ultra thick,inner sep=3pt,shape=ellipse,fit=(u1) (v1) (u2)] {};
    \end{pgfonlayer}
   \end{tikzpicture}
  \begin{tikzpicture}[
    vertex/.style={circle,draw=black,minimum size=24pt},
    edge/.style={very thick,color=black,double=green,double distance=30pt,line
    cap=round,draw opacity=0.2,fill opacity=1.0}]
    \node [vertex] (u1) {$u_1$};
    \node [vertex,draw=white,below=43pt of u1] {}; 
    \node [vertex,right=of u1] (v1) {$v$};
    \node [vertex,right=of v1] (u2) {$u_2$};
    \node [vertex,right=of u2] (v2) {$v'$};
    \node [vertex,right=of v2] (u4) {$u_4$};
    \node [vertex,right=of u4] (u3) {$u_3$};
    \node [vertex,right=of u3] (u5) {$u_5$};

    \begin{pgfonlayer}{background}
      \node [draw=black,ultra thick,inner sep=0pt,shape=ellipse,fit=(u1) (v1) (u2)] {};
      \node [draw=black,ultra thick,inner sep=0pt,shape=ellipse,fit=(u2) (v2) (u4)] {};
      \node [draw=black,ultra thick,inner sep=0pt,shape=ellipse,fit=(u4) (u3) (u5)] {};
    \end{pgfonlayer}
   \end{tikzpicture}
   }
   \caption{$\{u_1, u_2, u_3, u_4, u_5\}$ absorbs the pair $\{v, v'\}$.}\label{fig:absorber}
\end{figure}
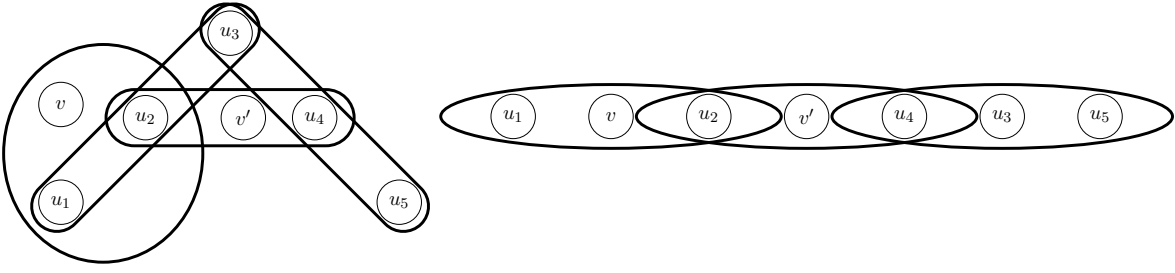
  \begin{lemma}{(Connecting Lemma)}\label{lem2}
   If $P = p_1 \dotsc p_{2a+1}$ and
   $Q = q_1 \dotsc q_{2b+1}$ are disjoint paths then
   $p_1 \dotsc p_{2a+1} v q_1 \dotsc q_{2b+1}$ is a path
   for every $v \in N(p_{2a+1}, q_1) \setminus (V(P) \cup V(Q))$
 \end{lemma}
 \noindent{\bf Proof.} This is clear from the definition of a loose path.
 \qed
 \begin{defn}
   A path $P$ with endpoints $v$ and $v'$ is said to absorb a
   set $W \subseteq V \setminus V(P)$ if there is a path
   $P'$ in $H$ with endpoints $v$ and $v'$ and such that $V(P')= V(P) \cup W$.
   A path $P$ is called $\gamma$-absorbing if $P$ absorbs every set $W
   \subseteq V\setminus V(P)$ such that $|W| \leq \gamma n$ and $|W|$ is even. 
 \end{defn}
 \begin{lemma}{(Absorbing Lemma)}\label{absorbing-lemma}
For every $0<\delta\leq 1/(24\cdot 10^3)$ there is $n_0$ such that for every
3-graph $H$ on $n> n_0$ vertices there is a path $P_{abs}$ in $H$ on at most $\delta n$ vertices which is $\delta^2$-absorbing. 
 \end{lemma}
 \noindent{\bf Proof.}  For  two distinct vertices $v, v' \in V$ let
 $\mathcal{A}(v,v')$ be the family of sets $S\in \binom{V}{5}$ that absorb
 $\{v,v'\}$. By Claim~\ref{lem1}, $|\mathcal{A}(v,v')| \geq \binom{n}{5}/500$.
 Let $\mathcal{F}$ be a family obtained by selecting every set from
 $\binom{V}{5}$ independently with probability $p := \delta n/(20
 \binom{n}{5})$. By the Chernoff bound \cite{chernoff}, with probability $1- o(1)$, the
 following two facts hold: (i) $|\mathcal{F}| \leq 2 p \binom{n}{5} = \delta
 n/10$; (ii) For every $\{v,v'\}$, $|\mathcal{A}(v,v') \cap \mathcal{F}| \geq
 p\binom{n}{5}/600 = \delta n/12000\geq 2 \delta^2 n$.
 
 The expected number of pairs $\{S_1, S_2\}$ such that $S_1, S_2 \in
 \mathcal{F}$ and $S_1 \cap S_2 \neq \emptyset$ is at most $5 \binom{n}{5}
 \binom{n}{4} p^2\leq 
 \delta^2 n/10$ and so by Markov's inequality, with probability at least
 $1/2$, the number of such pairs is at most $\delta^2 n/5$. Therefore, with
 positive probability, there exists a family $\mathcal{F}$ such that
 $|\mathcal{F}| \leq \delta n/10$, for every $\{v,v'\}$, $|\mathcal{A}(v,v')
 \cap \mathcal{F}| \geq 2 \delta^2 n$, and the number of $\{S_1, S_2\}$ such
 that $S_1, S_2 \in \mathcal{F}$ and $S_1 \cap S_2 \neq \emptyset$ is at most
 $\delta^2 n/5$. Let $\mathcal{F}'$ be obtained from $\mathcal{F}$ by
 deleting all intersecting sets and sets that do not absorb any $\{v,v'\}$.
 Then $|\mathcal{F}'| \leq \delta n/10$, and for every $\{v, v'\}$,
 $|\mathcal{A}(v,v') \cap \mathcal{F}'| > \delta^2 n$. 
 For each $S \in \mathcal{F'}$, $H[S]$ contains a path on five
 vertices, and using the co-degree condition and Lemma~\ref{lem2}
 we can connect the endpoints of these
 paths to obtain a new path $P_{abs}$. We have that $|V(P_{abs})| \leq  \delta
 n$.
 
 To show that $P_{abs}$ is $\delta^2$-absorbing, 
 consider $W \subseteq V \setminus V(P_{abs})$ such that $|W|$ is even and $|W| \le \delta^2 n$.
 Let $\mathcal{W} = \{\{w_1, w_1'\}, \dotsc, \{w_m,w_m'\}\}$ 
 be an arbitrary partition of $W$ 
 into sets of size $2$.
 We have that that $|\mathcal{A}(w_i,w_i') \cap \mathcal{F'}| > \delta^2 n$
 for every $i \in [m]$.
 Therefore, there is a matching between $\mathcal{W}$ and $\mathcal{F'}$
 so that every $\{w_i, w_i'\} \in \mathcal{W}$ is paired
 with some $S_i \in \mathcal{A}(w_i, w_i')$.
 This implies that $P_{abs}$ absorbs $W$ and the proof is complete.
\qed 
 
 Let $\delta := \min\{(\beta/1728)^2, 1/(24 \cdot 10^3)\}$, $\epsilon :=
 \delta/4$, and $C:=28/\delta^{5}$ 
 and let $P_{abs}$ be a path from Lemma~\ref{absorbing-lemma}.  Since $|V(P_{abs})| \leq \delta n$, for every two vertices $v, v' \in V$, 
 $|N(v,v') \cap (V\setminus V(P_{abs}))| \geq n/5$. Let $V_{r}$ be a set of
 size $\delta^{3} n$ selected uniformly at 
 random from $\binom{V \setminus V(P)}{\delta^{3} n}$.  In view of Chernoff's bound, we have
 \begin{equation}\label{resev}
 |N(v, v') \cap V_r| \geq 2 \delta^{4} n
 \end{equation}
 for every two vertices $v, v' \in V$.
 ($V_r$ will serve as the \emph{connecting reservoir}.)

 Our goal now is to construct a path $P$ in $H' := H[V \setminus V(P_{abs})]$
 such that $|V(H' - P) \setminus V_r| \le \delta^3 n$ and 
 $|V(P) \cap V_r| \le \delta^{4}n$. 
 This will prove Theorem~\ref{non-extremal}.
 Indeed, using \eqref{resev} and Lemma~\ref{lem2}, 
 we can connect $P_{abs}$ and $P$ 
 and then connect the endpoints of the resulting path to form a cycle $C$.
 Since $|V(C)| \ge (1 - \delta^2)n$,
 $P_{abs}$ can absorb $V(H - C)$ giving us a loose Hamilton cycle.
 
 Starting with $Q_0$ a longest path in $H' - V_r$,
 we will build $P$ iteratively, adding
 at least $\varepsilon C$ vertices in each step 
 with at most $7$ additional vertices from $V_r$.
 That is, we will construct paths $Q_0, \dotsc, Q_d$ in $H'$ so that
 $|V(H' - Q_d) \setminus V_r| \le \delta^3 n$,
 $|Q_i| \ge |Q_{i-1}| + \varepsilon C$ and
 $|V(Q_i - Q_{i-1}) \cap V_r| \le 7$ for $i \in [d]$.
 Note that we must have $d \le n/(\varepsilon C) \le \delta^4 n / 7$,
 so $|V(Q_i) \cap V_r| \le \delta^4 n$ for every $i \in [d]$.
 In particular, this shows that $Q_d$ will be the desired path $P$.

Let $Q$ be a path in $H'$ and let $W = V(H' - Q) \setminus V_r$.
We will now complete the proof by showing that if 
$|W| \ge \delta^3 n$ and $|V(Q) \cap V_r| \le \delta^4 n$ then
there exists a path $Q'$ in $H'$ such that
$|V(Q')| \geq |V(Q)| + \epsilon C$ and $|V(Q' - Q) \cap V_r| \le 7$.
%
To  find $Q'$ we partition $Q$ into blocks $B_1, \dots, B_l$ of consecutive
vertices so that the blocks are pairwise disjoint and for each $i$, 
$C+2 \leq |B_i| \leq (1+\epsilon) C$. Note that each $Q[B_i]$ contains a path of length at least $C$.
We first present a short claim and then consider three cases.
 \begin{claim}\label{lem3}
 Let $W$ be such that $|H[W]| \geq \gamma |W|^3$. Then there is a path on at least $\gamma |W|$ vertices in $H[W]$.
 \end{claim}
 \noindent{\bf Proof.}  Let $W' := W$.  If there is a vertex $w \in W'$ with
 $|N(w) \cap \binom{W'}{2}| < \gamma |W|^2$, then delete $w$ from $W'$.
 Continue if possible.  
 Note that the above process ends with $W' \neq \emptyset$ as otherwise,  $|H[W]| < \gamma |W|^2 \cdot |W|$. 
 So for every vertex $w \in W'$, $|N(w) \cap \binom{W'}{2}| \geq \gamma
 |W|^2$. Let $P$ be a longest path in $H[W']$. For an endpoint $v$ of $P$,
 note that if $e \in N(v) \cap \binom{W'}{2}$ then
 $e \cap V(P) \neq \emptyset$, because otherwise we could extend $P$.
 Therefore, if $|V(P)|< \gamma |W|$ then
 $|N(v) \cap \binom{W'}{2}| < \gamma |W|^2.$  \qed

 \bigskip
\noindent {\bf Case 1: $|H[W]| \geq \delta^{2}|W|^3$. }\\
By Claim~\ref{lem3} there is a path $P'$ in $H[W]$ on at least $
\delta^{2}|W|>C$ vertices and by Lemma~\ref{lem2}, we can connect one of the
endpoints of $Q$ with one of the endpoints of $P'$ using a vertex from
$V_r$ to form $Q'$.\\
{\bf Case 2}: There is a block $B_i$ such that $e(B_i, W, W)\geq (1/4 +
\epsilon)|B_i|\binom{|W|}{2}$.\\
A pair $\{w,w'\} \in \binom{W}{2}$ is called {\it good} if $|N(w,w')\cap B_i| \geq (1/4 + \epsilon/2)|B_i|$. Let $G_{good}$ be the graph induced by the good pairs. We have
$$|G_{good}| |B_i| + (1/4 + \epsilon/2)|B_i|\left(\binom{|W|}{2} -
|G_{good}|\right) \geq (1/4 + \epsilon) |B_i|  \binom{|W|}{2}$$
and so
$$|G_{good}| \geq \frac{\epsilon/2}{3/4 - \varepsilon/2}\binom{|W|}{2}= \Omega(|W|^2).$$
Color the edges of $G_{good}$ by setting $c(e)= N(e) \cap B_i$. 
Note that each $c(e)$ is a subset of $B_i$ of size
at least $(1/4 + \epsilon/2)|B_i|$ and since $|B_i| \leq 1.1 C$,  the number
of such subsets is less than $2^{1.1C}$. 
Therefore there exists $V_3 \subseteq B_i$ such that
if $F := c^{-1}(V_3)$ then $|F| \ge \Omega(|W|^2)$.
By the problem of
Zarankiewicz \cite{zaran}, the graph $F$ contains a complete bipartite
graph $K_{C, C}$. 
Let $V_1, V_2$ be the partite sets of $K_{C,C}$.
Then $H[V_1, V_2, V_3]$ is a complete $3$-partite $3$-graph. 
Since 
$|V_3| \geq (1/4 + \epsilon/2)|B_i|$ and $|V_1|, |V_2| \geq C$, 
$H[V_1, V_2, V_3]$
contains a path $P'$ on at least $(1+ 2\epsilon)|B_i|$ vertices. We connect
the endpoints of $P'$ with two paths obtained by removing $B_i$ from $Q$ using
$2$ vertices from $V_r$.  \\
{\bf Case 3}: For every block $B_i$,  $e(B_i, W, W)\leq (1/4 +
\epsilon)|B_i|\binom{|W|}{2}$.\\
Since we are not in Case 1,
\begin{equation*}
  e(V(Q), W, W) \geq (n/4 - |V_r| - |P_{abs}|)\binom{|W|}{2} - \delta^2|W|^3
  \geq (1/4 - 2\delta) n \binom{|W|}{2}.
\end{equation*}
A block $B$ is called {\it good} if
$e(B, W, W) \geq (1/4 -2\sqrt{\delta})|B|\binom{|W|}{2}$.  

Let $l'$ be the number of good blocks. If $l' < (1-3\sqrt{\delta})n/C$ then
\begin{equation*}
  \begin{split}
    e(V(Q), W, W) &\leq l' \cdot (1/4 + \epsilon)(1 + \varepsilon)C\binom{|W|}{2} +  
    (l-l')(1/4 - 2\sqrt{\delta})(1 + \varepsilon)C \binom{|W|}{2} \\
    &\leq l' \cdot (1/4 + 2\epsilon)C\binom{|W|}{2} +  
    (l-l')(1/4 - \sqrt{\delta})C \binom{|W|}{2} \\
   &< (2 \epsilon + \sqrt{\delta})(1 - 3 \sqrt{\delta})n\binom{|W|}{2} + (1/4
   - \sqrt{\delta}) lC\binom{|W|}{2} 
   < (1/4 - 2 \delta) n \binom{|W|}{2}
 \end{split}
\end{equation*}
which is not possible. Thus at least $(1 -3 \sqrt{\delta})n/C$ blocks are good.
Using the same argument as in the previous case for every good block $B_i$ we
can find $V_1(i), V_2(i) \subset W$, $V_3(i) \subset B_i$ such that
$|V_3(i)|= \left\lceil(1/4 - 3\sqrt{\delta}) C\right\rceil$, $|V_1(i)|,
|V_2(i)| \ge 9C/4$, and
$H[V_1(i), V_2(i), V_3(i)]$ is a complete $3$-partite $3$-graph. Let
$U:= \bigcup (B_i\setminus V_3(i))$ where the union is taken over all $i$ such that $B_i$ is good. Then
$|U| \geq 3n/4$ and, since we are not in the extremal case,
$|H[U]| \geq \beta n^3$. Thus there exist three good blocks, say $B_1, B_2,
B_3$ and $Z \subseteq B_1 \cup B_2 \cup B_3$ such that $|U \cap Z| = 9C/4$
and $|H[U \cap Z]| \geq (\beta/2) \cdot
(3C/4)^3 = (\beta/54) \cdot (9C/4)^3$. Now, by Lemma~\ref{lem3}, there is a path on 
$(\beta/54)(9C/4) = \beta C/24$ vertices in
$H[U \cap Z]$ and, since these blocks are good, there are
three vertex disjoint paths in $H[(Z \setminus U) \cup W]$, 
each on at least $(1 -12 \sqrt{\delta})C$ vertices. Altogether the paths have
at least $3C + \beta C/48 >  (1+ \beta/150)(|B_1| + |B_2|+|B_3|)$ vertices. We can form $Q'$ by
connecting the four new paths with $V(Q) \setminus (B_1 \cup B_2 \cup B_3)$ 
using at most seven vertices from $V_r$. \qed
\end{section}
\section{Extremal case}
In this section, we prove the extremal case.
\begin{theorem}\label{extremal}
There is $\beta_0>0$ such that for every $0<\beta< \beta_0$ there is an
integer $n_0$ such that every 3-graph $H = (V, E)$ on $n\geq n_0$ vertices
where $n$ is even which is $\beta$-extremal and with minimum co-degree at least $n/4$ contains a loose Hamilton cycle.
\end{theorem}
\noindent We first  present a rough outline of the proof.  
The following definition is critical. 
\begin{defn}
  Let $H = (V, E)$ be a $3$-graph and $V^* \subseteq V$.
  A vertex $v \in V$ is
  \emph{$\gamma$-good for $V^*$} if 
  \begin{equation*}
    \left|N(v) \cap \binom{V^*}{2}\right| \ge \left(1 - \gamma\right)\binom{|V^*|}{2}
  \end{equation*}
  and a pair of vertices $\{v, v'\} \in \binom{V}{2}$ is 
  \emph{$\gamma$-good for $V^*$} if 
  \begin{equation*}
    \left|N(v, v') \cap V^*\right| \ge \left(1 - \gamma\right)|V^*|.
  \end{equation*}
\end{defn}
For some small $\beta$ we are given a $\beta$-extremal $3$-graph $H = (V,E)$
on $n$ vertices with minimum co-degree $n/4$.  So there exists 
a partition $\{U, W\}$ of $V$ such that $|W| = \left\lceil 3n/4 \right\rceil$
and $H[W]$ has at most $\beta |V|^3$ edges.
Say a vertex is \emph{$\gamma$-exceptional} if it is in $U$ and is not
$\gamma$-good for $W$ or if it is in $W$ and is in less than $(1-\gamma)|W|$ of
the pairs in $\binom{W}{2}$ that are $\gamma$-good for $U$.
If we choose $\gamma$ to be a positive constant that is small relative 
to $\beta$ then only a small fraction of the vertices will be $\gamma$-exceptional.
We first construct a loose path, which we will call $P_s$, that contains
all of these $\gamma$-exceptional vertices.
We also require $|W \setminus V(P_s)| = 3(|U \setminus V(P_s)| - 1)$,
$|P_s|$ to be a small fraction of $n$
and the endpoints of $P_s$ to be vertices in $U$ that are $\gamma$-good for $W$.
We then remove all of $P_s$ except the endpoints from $H$ to create
the graph $H'$.
Note that $H'$ is very well behaved in the sense that
there are no $\gamma$-exceptional vertices in $H'$ 
and $3(|U \cap V(H')| - 1) = |W \cap V(H')|$. 
We can therefore argue that there exists a spanning loose path of $H'$
that has the same endpoints as $P_s$ and in which every edge
has exactly one vertex in $U$ and two vertices in $W$.
Combining $P_s$ and the spanning path of $H'$ gives a Hamilton cycle of $H$.

We begin the formal proof with a general lemma that is used to produce the
spanning path of $H'$ in the preceding description.
It relies on the following theorem of K\"{u}hn and Osthus \cite{KO-match}.
\begin{theorem}\label{matching}
For all positive constants $d, v_0 , \eta \leq 1$ there is a positive $\epsilon =
\epsilon(d, v_0 , \eta)$ and an integer $N_0 = N_0(d, v_0 , \eta)$ such that the following holds for all $n \geq N_0$
and all $v\geq v_0$ . Let $G = (A, B)$ be a $(d, \epsilon)$-superregular bipartite graph whose vertex
classes both have size $n$ and let $F$ be a subgraph of $G$ with $e(F) = ve(G)$. Choose a
perfect matching $M$ uniformly at random in $G$. Then with probability at
least $1 -e^{-\epsilon n}$
we have
$$(1 - \eta)vn \leq |M \cap E(F)| \leq (1 + \eta)vn.$$
\end{theorem}

\begin{lemma}\label{main-extremal}
There exists $\gamma>0$ and $m_0$ such that the following holds. 
Let $H' = (V,E)$ be a $3$-graph.
If there exists $\{U^*, W^*\}$ a partition of $V$ such that 
\begin{itemize}
  \item $|U^*| = m$ and $|W^*| = 3(m - 1)$ for some $m \ge m_0$;
  \item every vertex $u\in U^*$ is $\gamma$-good for $W^*$; and
  \item for every vertex $w_1\in W^*$ all but at most $\gamma |W^*|$ vertices 
    $w_2 \in W^*$ are such that $\{w_1, w_2\}$ is $\gamma$-good for $U^*$.
\end{itemize}
Then for any $u, u' \in U^*$ there exists a spanning loose path
of $H'$ with endpoints $u$ and $u'$. 
\end{lemma}
\noindent Let $d:=1$, $\nu_0 := 19/20$ and $0 < \eta \le 3/40$.
  There exists $\epsilon = \epsilon(d, \nu_0, \eta)$ such
  that Theorem~\ref{matching} holds. 
  Let $\gamma := \min(\epsilon^2/3, 1/100)$.
  Let $G$ be a graph on $W^*$ such that
  $E(G) := \{ww' : \{w,w'\} \text{ is $\gamma$-good for $U^*$ } \}$.
 Partition $W^*$ into three sets $W_1, W_2, W_3$, each of order $m-1$.
  Since $\delta(G) \ge (1 - \gamma)|W^*| \ge (3 - \varepsilon^2)(m-1)$,
  $G[W_1, W_2]$ and $G[W_2, W_3]$ are both 
  $(d, \epsilon)$-superregular bipartite graphs.
  For $u \in U^*$, let $F^1_u:=\{ww' \in G[W_1, W_2]: \text{$\{u,w,w'\} \in
H$}\}$ and let
$F^2_u:=\{ww' \in G[W_2, W_3]: \text{$\{u,w,w'\} \in H$}\}$. Note that $|F^1_u| > (1- 5\gamma)|G[W_1, W_2]|$ and $|F^2_u| > (1- 5\gamma)|G[W_2, W_3]|$. Therefore, by Theorem~\ref{matching}, there is a matching $M_1$ in $G[W_1, W_2]$ and a matching $M_2$ in $G[W_2,W_3]$,
  such that every vertex in $U^*$ is in a hyperedge with all but 
  less than $(5 \eta + \gamma)(m-1)$ of the edges in $M_1$ and all but 
  less than $(5 \eta + \gamma)(m-1)$ of the edges in $M_2$.
  Label the vertices of $W^*$ so that
  $W_1 = \{a_1, \dotsc, a_{m-1}\}$,
  $W_2 = \{b_1, \dotsc, b_{m-1}\}$ and 
  $W_3 = \{c_1, \dotsc, c_{m-1}\}$; 
  and so that 
  $M_1 = \{a_1b_1, \dotsc, a_{m-1}b_{m-1}\}$ and 
  $M_2 = \{b_1c_1, \dotsc, b_{m-1}c_{m-1}\}$.
  Construct an auxiliary bipartite graph $L$ with one part 
  $U^*$ and the other the set of length two paths
  $\{a_1b_1c_1, \dotsc, a_{m-1}b_{m-1}c_{m-1}\}$
  in which an edge is present between
  $u \in U^*$ and $a_ib_ic_i$ if and only if
  $\{u, a_i, b_i\}$ and $\{u, b_i, c_i\}$  are both in the hypergraph $H$.
  For every $i \in [m - 1]$,
  $\{a_i,b_i\}$ and $\{b_i,c_i\}$ are both $\gamma$-good
  for $U$, so $a_ib_ic_i$ has degree at least $(1 - 2\gamma)m$ in $L$.
  Furthermore, by construction every $u \in U$ has degree at least 
  $(1 - (10 \gamma + 2 \eta))(m-1) \ge (3/4)(m-1)$. This
  implies that for any pair $u, u' \in U^*$.
  there exists a spanning path of $L$ with endpoints $u$ and $u'$.
  This path corresponds to a spanning loose path of the hypergraph $H'$ 
  with endpoints $u$ and $u'$. \qed


We now present the proof of the extremal case.

\noindent{\bf Proof of Theorem~\ref{extremal}.} 
The value of $\beta_0$ follows from the computations
and the application of Theorem~\ref{matching}. In what follows we will assume
that $\beta>0$ is sufficiently small and will let $\sigma := (50 \beta)^{1/4}$. 

Because $H$ is $\beta$-extremal, there is a set $W \subset V$ such that $|W| \in \{3n/4,
(3n+2)/4\}$ and $|H[W]| \leq \beta n^3$. Let $U:= V \setminus W$. 
We have $e(U, W,W)\geq \frac{n}{4} \binom{|W|}{2} - 3 \beta n^3\geq
(1 - \sigma^4)|U|\binom{|W|}{2}$.
Consequently all but at most $\sigma^2|U|$ vertices $u\in U$ are
$\sigma^2$-good for $W$. Indeed, if $U'$ is the set of $\sigma^2$-good vertices, then 
\begin{multline*}
  (1-\sigma^4)|U|\binom{|W|}{2} \leq |U'| \binom{|W|}{2} + 
  (|U|- |U'|)(1-\sigma^2)\binom{|W|}{2} 
  = ((1-\sigma^2)|U| + \sigma^2|U'|) \binom{|W|}{2}
\end{multline*}
and so $|U'| \geq (1- \sigma^2)|U|$.

Let $G_{good}$ be a graph on $W$ induced by the pairs that are $\sigma^2$-good. We have
$$(1-\sigma^4)|U|\binom{|W|}{2} \leq |G_{good}||U| + \left(\binom{|W|}{2} -
|G_{good}|\right) \left(1- \sigma^2\right)|U|$$
and so $|G_{good}| \geq (1- \sigma^2)\binom{|W|}{2}$. Thus all but at
most $\sigma |W|$ vertices in $w\in W$ have $\deg_{G_{good}}(w) \geq
(1- \sigma)|W|$.

A vertex $v\in V$ is called {\it $\alpha$-acceptable for $W$} if for at least $\alpha n$ vertices $w \in W$, $|N(v,w)\cap U| \geq \alpha n$.
A vertex $v\in V$ is called {\it $\alpha$-acceptable for $U$} if $|N(v) \cap
\binom{|W|}{2}| \geq \alpha n^2$.
\begin{claim}
Every vertex $v\in V$ is either $1/16$-acceptable for $W$ or is $1/16$-acceptable for $U$.
\end{claim}
 \noindent{\bf Proof.} 
 We have that $\sum_{w \in W \setminus \{v\}} |N(v,w)| \geq (|W| - 1)n/4$, and if $|N(v) \cap \binom{W}{2}| \leq \alpha n^2$,
 then $e(\{v\}, W, U) \geq (3/16 -2 \alpha)n^2$. Suppose that for less than $\alpha n$ vertices $w$ in $W$, 
$|N(v,w) \cap U| \geq \alpha n$. 
Then
\avoidbreak
$$(3/16 -2\alpha) n^2 < \alpha n |U| + (|W| - \alpha n) \alpha n < \alpha n^2,$$
so $\alpha > 1/16$. \qed

To continue with the proof, let $\widetilde{U}\subset U$ be the set of vertices that
are not $\sigma^2$-good for $W$, and let $\widetilde{W}$ be the set of vertices $w
\in W$ such that $\deg_{G_{good}}(w) < (1- \sigma)|W|$. Note that
$|\widetilde{U}\cup \widetilde{W}| \leq \sigma n$ and distribute vertices from 
$\widetilde{U}\cup \widetilde{W}$ to $U', W'$ so that a vertex assigned to $U'$ (resp. $W'$) is
$1/16$-acceptable for $U$ (resp. $W$).  
If $v\in \widetilde{U}\cup \widetilde{W}$ is
$1/16$-acceptable for both $U$ and $W$ then add $v$ to $U'$.
Let $U'':= U \setminus \widetilde{U}$ and $W'' := W \setminus \widetilde{W}$. 
Note that $(1/4 - 2\sigma)n < |U' \cup U''| < (1/4 + 2\sigma)n$.

If $|U' \cup U''| < \frac{n}{4}$ then let 
$p_1 := \left\lceil \frac{n}{4} \right\rceil - |U' \cup U''|$
else let
$p_1 := 0$.
Define $p_2 := 2 \left( |U' \cup U''| + p_1 - \frac{n}{4}\right)$.
Note that $p_1$ and $p_2$ are integers,
$0 \le p_1 \le 2 \sigma n$ and $0 \le p_2 \le 4 \sigma n$
and that
\begin{equation*}
  3\left(|U' \cup U''| + p_1 \right) - 2 p_2 =
  n - |U' \cup U''| - p_1 = |W' \cup W''| - p_1,
\end{equation*}
so $3|U' \cup U''| + 4p_1 -2p_2 = |W' \cup W''|$.
Suppose then that we have
constructed a path $P_s$ with
endpoints $u, u' \in U''$ such that
\begin{equation*}
  3(|V(P_s - u - u') \cap  (U' \cup  U'')| - 1) + 4p_1 - 2p_2 = 
  |V(P_s) \cap  (W' \cup  W'')|,
\end{equation*}
$U'\cup W' \subset V(P_s)$ and $|P_s| \le 10 \sigma n$.
If we let $H' := H[\left(V \setminus V(P_s)\right) \cup \{u, u'\}]$
and let $m := |V(H') \cap (U' \cup U'')|$ then 
\begin{equation*}
  3(m - 1) = 3|U' \cup U''| - 3(|(V(P_s) - u - u') \cap (U' \cup U'')| - 1) =
  |V(H') \cap (W' \cup W'')|.
\end{equation*}
Therefore we can  apply Lemma~\ref{main-extremal}, to $H'$ 
provided $\sigma \le \frac{\gamma}{50}$.
This gives a Hamilton cycle of $H$ thereby proving Theorem~\ref{extremal}.

\begin{figure}
  \centering
  \scalebox{1.0} {
  \begin{tikzpicture}[
    vertex/.style={circle,fill=black,minimum size=5pt,inner sep=0pt}
    ]
    \node [vertex] (v1) {};
    \node [vertex,below right=8pt of v1] (v2) {};
    \node [vertex,below right=8pt of v2] (v3) {};
    \node [rotate fit=-45,draw=black,inner sep=5pt,
    shape=rounded rectangle,fit=(v1) (v2) (v3)] {};
    \UpEdge{3}
    \TwoPath{5}
    \DownEdge{9}
    \UpEdge{11}
    \DownEdge{13}
    \UpEdge{15}
    \TwoPath{17}
    \DownEdge{21}
    \UpEdge{23}

    \draw [dashed,very thick] ($(v3) + (-75pt,5pt)$) -- ($(v23) + (50pt,5pt)$);
    \node (Q1) [inner sep=5pt,thick,dashed,draw=black,shape=rounded rectangle,
    fit=(v5) (v6) (v7) (v8) (v9)] {};
    \node (Q2) [inner sep=5pt,thick,dashed,draw=black,shape=rounded
      rectangle, fit=(v17) (v18) (v19) (v20) (v21)] {};
    \node (Q1) at ($(Q1) + (0pt,15pt)$) {$Q_1$};
    \node (Q2) at ($(Q2) + (0pt,15pt)$) {$Q_2$};

    \node at ($(v3) + (-55pt, 15pt)$) {$W' \cup W''$};
    \node at ($(v3) + (-55pt, -3pt)$) {$U''$};
  \end{tikzpicture}
  }
  \caption{The path $P_1$ when $p_1 = 2$, $u_i = 4$ and $w_i = 21$.
  Every vertex that is not in some $Q_i$ 
  is in $U'' \cup W''$\label{fig:path_P1}}
\end{figure}
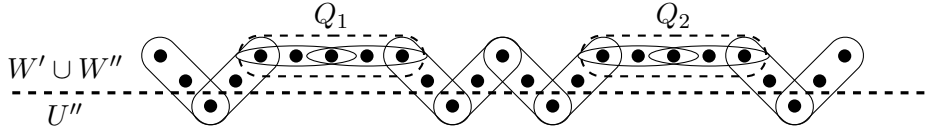

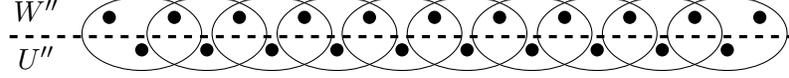
\begin{figure}
  \centering
  \scalebox{1.0} {
  \begin{tikzpicture}[
    vertex/.style={circle,fill=black,minimum size=5pt,inner sep=0pt}
    ]
    \node [vertex] (v1) {};
    \node [vertex,below right=12pt of v1] (v2) {};
    \node [vertex,above right=12pt of v2] (v3) {};
    \node [draw=black,inner sep=0pt, 
    shape=ellipse,inner sep=1pt,fit=(v1) (v2) (v3)] {};
    \UpDownEdge{3}
    \UpDownEdge{5}
    \UpDownEdge{7}
    \UpDownEdge{9}
    \UpDownEdge{11}
    \UpDownEdge{13}
    \UpDownEdge{15}
    \UpDownEdge{17}
    \UpDownEdge{19}
    
    \draw [dashed,very thick] ($(v2) + (-50pt,5pt)$) -- ($(v18) + (50pt,5pt)$);
    \node at ($(v2) + (-40pt, 15pt)$) {$W''$};
    \node at ($(v2) + (-40pt, -3pt)$) {$U''$};
  \end{tikzpicture}
  }
  \caption{The path $P_2$ when $p_2 = 10$, $u_2 = 10$ and $w_2 = 11$.}\label{fig:path_P2}
\end{figure} 

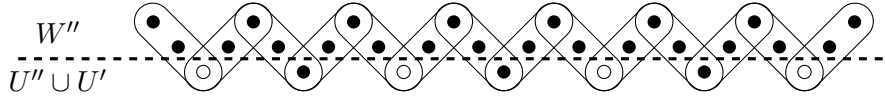
\begin{figure}
  \centering
  \scalebox{1.0} {
  \begin{tikzpicture}[
    vertex/.style={circle,fill=black,minimum size=5pt,inner sep=0pt}
    ]
    \node [vertex] (v1) {};
    \node [vertex,below right=8pt of v1] (v2) {};
    \node [vertex,fill=white,draw=black,below right=8pt of v2] (v3) {};
    \node [rotate fit=-45,draw=black,inner sep=0pt, 
    shape=rounded rectangle,inner sep=5pt,fit=(v1) (v2) (v3)] {};
    \UpEdge{3}
    \DownEdge{5}
    \UpEdge{7}
    \DownEdgeWhite{9}
    \UpEdge{11}
    \DownEdge{13}
    \UpEdge{15}
    \DownEdgeWhite{17}
    \UpEdge{19}
    \DownEdge{21}
    \UpEdge{23}
    \DownEdgeWhite{25}
    \UpEdge{27}
    
    \draw [dashed,very thick] ($(v3) + (-70pt,5pt)$) -- ($(v27) + (30pt,5pt)$);
    \node at ($(v3) + (-55pt, 15pt)$) {$W''$};
    \node at ($(v3) + (-55pt, -3pt)$) {$U'' \cup U'$};
  \end{tikzpicture}
  }
  \caption{The path $P_3$ when $|U'| = 4$, $u_3 = 7$ and $w_3 = 22$.  The
  white vertices are the  vertices of $U'$.}\label{fig:path_P3}
\end{figure} 

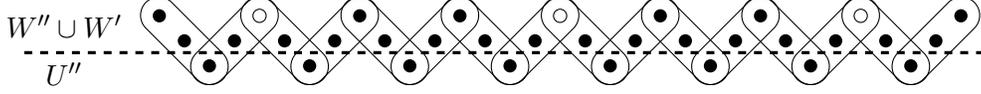
\begin{figure}
  \centering
  \scalebox{1.0} {
  \begin{tikzpicture}[
    vertex/.style={circle,fill=black,minimum size=5pt,inner sep=0pt}
    ]
    \node [vertex] (v1) {};
    \node [vertex,below right=8pt of v1] (v2) {};
    \node [vertex,below right=8pt of v2] (v3) {};
    \node [rotate fit=-45,draw=black,inner sep=0pt, 
    shape=rounded rectangle,inner sep=5pt,fit=(v1) (v2) (v3)] {};
    \UpEdgeWhite{3}
    \DownEdge{5}
    \UpEdge{7}
    \DownEdge{9}
    \UpEdge{11}
    \DownEdge{13}
    \UpEdgeWhite{15}
    \DownEdge{17}
    \UpEdge{19}
    \DownEdge{21}
    \UpEdge{23}
    \DownEdge{25}
    \UpEdgeWhite{27}
    \DownEdge{29}
    \UpEdge{31}
    
    \draw [dashed,very thick] ($(v3) + (-70pt,5pt)$) -- ($(v31) + (30pt,5pt)$);
    \node at ($(v3) + (-55pt, 15pt)$) {$W'' \cup W'$};
    \node at ($(v3) + (-55pt, -3pt)$) {$U''$};
  \end{tikzpicture}
  }
  \caption{The path $P_4$ when $|W' \setminus V(P_1)| = 3$, $u_4 = 8$ and $w_4 = 25$.
  The white vertices are the vertices of $W' \setminus V(P_1)$.}\label{fig:path_P4}
\end{figure}

We first give a brief overview of the construction of such a $P_s$ before proceeding
with the details.
First at most four disjoint paths each with
their endpoints in $W''$ are constructed: $P_1$, $P_2$, $P_3$ and $P_4$.
Define $u_i := |V(P_i) \cap (U' \cup U'')|$
and $w_i := |V(P_i) \cap (W' \cup W'')|$.
If $p_1 > 0$, then $|W' \cup W''| > \left\lfloor 3n/4 \right\rfloor$.
In this case, we will first find $p_1$ disjoint paths
of length two, $Q_1, \dotsc, Q_{p_1}$, in $H[W' \cup W'']$
and then connect these paths to form $P_1$ (see Figure~\ref{fig:path_P1}).
$P_1$ will be constructed so that 
$V(P_1) \subset U'' \cup W' \cup W''$,
$u_1 = 2p_1$
and $w_1 = 3u_1 + 1 + 4p_1$.
If $p_2 > 0$, then either $|U' \cup U''| > \left\lceil{n/4}\right\rceil$ or $n \equiv 2 \pmod 4$. 
In this case, $P_2$ will be constructed so that 
$V(P_2) \subset U'' \cup W''$
and such that $u_2 = p_2$
and $w_2 = 3u_2 + 1 - 2p_2$ (see Figure~\ref{fig:path_P2}).
If $U'\neq \emptyset$,
$P_3$ will be constructed so that 
$U' \subset V(P_3)$,
$u_3 = 2 |U'|-1$ 
and $w_3 = 3u_3 + 1$ (see Figure~\ref{fig:path_P3}).
Finally, if $W' \setminus V(P_1) \neq \emptyset$,
$P_4$ will be constructed so that 
$W' \setminus V(P_1) \subset V(P_4)$,
$u_4 = 3 |W' \setminus V(P_1)|-1$ and $w_4 = 3u_4 + 1$ (see Figure~\ref{fig:path_P4}).  

We will then connect all of the paths 
$P_1$, $P_2$, $P_3$ and $P_4$ which have been 
constructed to form $P_s$.
For $i \in [3]$. let  $x, y \in W''$ be endpoints of 
$P_i$ and $P_{i+1}$ respectively.
We join $P_i$ with $P_{i+1}$ 
by selecting unused vertices $x', y' \in W''$ such that
$|N(x, x') \cap N(y, y') \cap U| \ge (1 - 2 \sigma)|U|$
so we can greedily select an unused $z \in U''$
such that $\{x, x', z\}, \{z, y', y\} \in H$.
Therefore, we will use two unused vertices from $W''$ and
one unused vertex from $U''$ at each juncture.
We can then greedily add one unused vertex
from $W''$ and one unused vertex from $U''$ to each end to create
the desired path $P_s$.  

When $U' \cup W' = \emptyset$,
$|U''| = \frac{n}{4}$ and $|W''| = \frac{3n}{4}$
none of the four paths are not constructed. 
In this case, we let 
$P_s := \{\{u_1, w_1, w_2\}, \{w_2, w_3, u_2\}\}$
for some $u_1, u_2 \in U''$ and $w_1, w_2, w_3 \in W''$.

Now we complete the proof by showing that we can construct the four
paths $P_1$, $P_2$, $P_3$ and $P_4$ as described above.

So suppose $p_1 > 0$ and note that 
because $|U'\cup U''| = \left\lceil \frac{n}{4} \right\rceil - p_1$, 
$H[W' \cup W'']$ has minimum co-degree $p_1$.
\begin{lemma}\label{indep-paths}
There are at least $p_1$ vertex-disjoint paths of length two in 
$H[W' \cup W'']$.
\end{lemma}
 \noindent{\bf Proof.} 
 For any $Z \subseteq V$ if
 $H[Z]$ does not contain a path with two edges, then it is
 not hard to see that every nontrivial
component of $H[Z]$ either has less than $5$ vertices or
is a double star. A double star
is a hypergraph $S$ where $S = \{\{u,w, z_i\}| i=1, \dots k\}$.  
This implies that $e(H[Z]) \le |Z|$. 

 Let $\mathcal{M}$ be a collection of vertex-disjoint
 paths of length two of the maximum size and suppose that $|\mathcal{M}| <
 p_1$.  Since $e(H[W' \cup W'']) \ge p_1\binom{|W' \cup W''|}{2}/3 > |W' \cup W''|$, 
 we can assume $|\mathcal{M}| \ge 1$.
 Let $Z' := \bigcup_{P \in \mathcal{M}} V(P)$ and $Z:= \left(W'\cup
 W''\right) \setminus  Z'$.
If there is a path $P \in \mathcal{M}$ such that for two vertices $v_1, v_2
\in V(P)$, $|N(v_1) \cap \binom{Z}{2}| \geq 5|Z|$ and $|N(v_2) \cap
\binom{Z}{2}| \geq 5|Z|$, then we can increase $|\mathcal{M}|$. Indeed,
$N(v_1) \cap \binom{Z}{2}$ contains a matching $M_1$ of size two and there are
less than $4|Z|-4$ pairs in $N(v_2) \cap \binom{Z}{2}$ incident to vertices from
$V(M_1)$. Thus $N(v_2) \cap \binom{Z}{2}$ contains a matching of size two that does not contain vertices of $M_1$.

Therefore, for every $P \in \mathcal{M}$, $|N(V(P)) \cap \binom{Z}{2}| <
\binom{|Z|}{2} + 20 |Z|$ and, consequently, 
$e\left(Z', Z, Z\right) < 
|\mathcal{M}|\left( \binom{|Z|}{2} + 20 |Z|\right)$.
Furthermore, 
$$e(H[Z])
> \frac{1}{3}\left(p_1 \binom{|Z|}{2} -  
|\mathcal{M}|\left( \binom{|Z|}{2} + 20 |Z|\right)\right) 
\geq \frac{1}{10} \binom{|Z|}{2}$$
as $|\mathcal{M}| \le p_1 -1/2$ and $p_1 \le 2 \sigma n$. 
Therefore $e(H(Z)) > |Z|$ so there exists a path of length $2$ in $H[Z]$ which
contradicts the maximality of $\mathcal{M}$.
\qed

Let $Q_1, Q_2, \dots, Q_{p_1}$ be the paths from Lemma~\ref{indep-paths}.
We connect the $Q_i$'s to form $P_1$ as follows. 
Suppose that $Q_1, \dots, Q_i$ have been connected into a path and $v_i$ is
the endpoint of the constructed path  and the endpoint of  $Q_i$. Since $v_i$
is $1/16$-acceptable for $W$ we can select $v, w, w' \in W''$ and $u \in U''$
not previously used such that $\{v_i, w', u\}, \{u,w, v\}\in H$. Let $v_{i+1}$
be one of the endpoints of $Q_{i+1}$. Since $v \in W''$ and 
$v_{i+1}$ is $1/16$-acceptable for $W$ we can find unused 
$x, x' \in W''$ and $u' \in
U''$ and such that  $\{v, x, u'\}, \{u', x', v_{i+1}\} \in H$.
By a similar argument, we can add 
three unused vertices from $W''$ and an unused vertex from $U''$ to
both ends of the resulting path so that the endpoints of $P_1$ are in $W''$. 

If $p_2 > 0$, let $P_2 =
\{v_1, u_1, v_2\}, 
\{v_2, u_2, v_3\} \dots, 
\{v_{p_2}, u_{p_2}, v_{{p_2}+1}\}$ be a path on $H - P_1$ such that 
$u_1, \dotsc, u_{p_2} \in U''$ and 
$v_1, \dotsc, v_{p_2 + 1} \in W''$. 
The path $P_2$ can be easily constructed as for 
$v_i$ we can select an unused $v_{i+1}$ with 
$|N(v_i, v_{i+1})\cap U| \geq (1- 2 \sigma)|U|$.

Build paths $P(u)$ for every $u\in U'$ one by one as follows. For the general
step, consider $u\in U'$ and find two edges $\{v_1, v_2,  u\}$, $\{u,v_3, v_4\}$
containing unused vertices $v_1, v_2, v_3, v_4 \in W''$. Note that this is
possible as $|U'| \leq 2 \sigma n$ and every vertex 
$u\in U'$ has $|N(u) \cap \binom{W''}{2}| \geq n^2/20$.
 Since  the endpoints of $P(u)$ are in $W''$, we can connect  $P(u)$ with
 $P(u')$  by using two vertices from $W''$ and one from $U''$. Let $P_3$ be
 the path obtained by connecting all of the $P(u)$'s.
 Construct $P_4$ in the similar fashion by considering paths $P(w)$ consisting of four edges 
 $\{v_1,v_2, u_1\}$, $\{u_1,v_3, w\}$, 
 $\{w, v_4, u_2\}$, $\{u_2, v_5, v_6\}$ with 
$u_1, u_2 \in U''$ and $v_1, v_2, v_3, v_4, v_5, v_6\in W''$ all unused.

\section{Acknowledgements}
We would like to thank the anonymous referees for their helpful comments.

\hbadness 10000\relax
\providecommand{\bysame}{\leavevmode\hbox to3em{\hrulefill}\thinspace}
\providecommand{\MR}{\relax\ifhmode\unskip\space\fi MR }
\providecommand{\MRhref}[2]{%
  \href{http://www.ams.org/mathscinet-getitem?mr=#1}{#2}
}
\providecommand{\href}[2]{#2}



\end{document}